\documentclass[12pt,a4paper]{amsart}
\pdfoutput=1
\usepackage{latexsym}
\usepackage{amsthm}
\usepackage{amsmath}
\usepackage{amssymb}
\usepackage{ulem}
\usepackage{color}
\usepackage{pgfplots}
\pgfplotsset{compat=1.5}
\usepackage{hyperref}
\advance\hsize by 14pt
\newskip\abstractindent         \abstractindent=3pc
\long\def\block #1\endblock{\vskip 6pt
        {\leftskip=\abstractindent \rightskip=\abstractindent
        \noindent #1\endgraf}\vskip 6pt}
 
\long\def\ext #1\endext{\removelastskip\block #1\endblock}

\def\today{\ifcase\month\or January \or February\or
March\or April\or May\or June\or July\or August\or
September\or October\or November\or December\fi
\space\number\day, \number\year}
 \newdimen\fullhsize
\fullhsize=\hsize
\def\fullline{\hbox to\fullhsize}

\def\vol{{\text{\rm vol}}}

\chardef\other=12

\newtheorem{theorem}{Theorem} 

\newtheorem{conjecture}[theorem]{Conjecture}

\numberwithin{theorem}{section} \numberwithin{equation}{section}

\begin{document}
\title{Aliquot sequences\\ with small starting values}
\author{Wieb Bosma}
\email{bosma@math.ru.nl}
\address{Radboud Universiteit, Heijendaalseweg 135, 6525 AJ Nijmegen, The Netherlands}

\begin{abstract}
We describe the results of the computation of aliquot sequences
with small starting values. In particular all sequences with
starting values less than a million have been computed until either
termination occurred (at 1 or a cycle), or an entry of 100 decimal 
digits was encountered. All dependencies were recorded, and numerous 
statistics, curiosities, and records are reported.
\end{abstract}
\maketitle
\section{Introduction}
\noindent
{\it Aliquot sequences} arise from iterating the sum-of-proper-divisors function
$$s(n)=\sum_{\genfrac{}{}{0pt}{}{d\vert n}{d<n}}d,$$
assigning to an integer $n>1$ the sum of its {\it aliquot} divisors
(that is, excluding $n$ itself). 
Iteration is denoted exponentially, so $s^k$ is shorthand for
applying $k\geq1$ times the function $s$.
We say that an aliquot
sequence {\it terminates (at 1)} if $s^k(n)=1$ for some $k$; this happens when
and only when $s^{k-1}(n)$ is prime. It is possible that $s^{k+c}(n)=s^k(n)$,
for some $c>0$ and all $k\geq k_0$, that is,
to hit an {\it aliquot cycle of length $c$}, where we take $c>0$ 
and $k_0$ minimal.
Case $c=1$ occurs when $n$ is a perfect number (like 6),
and $c=2$ when $n\neq m$ 
form a pair of {\it amicable numbers}:
$s(n)=m$ and $s(m)=n$.
See Section \ref{sec:cycles} for more cycles.

The main open problem regarding aliquot sequences is the 
conjecture attributed to Catalan \cite{Catalan} and Dickson \cite{Dickson}.
\begin{conjecture}
All aliquot sequences remain bounded.
\end{conjecture}
\medskip\noindent
If true, it would imply that for every $n$ after finitely many steps we either
hit a prime number (and then terminate at 1) or we find an aliquot cycle.
Elsewhere we comment upon some of the heuristics to support or refute
this conjecture \cite{BK}. 

We will call an aliquot sequence
{\it open} if it is not known to remain bounded.
This notion depends on our state of knowledge. The point of
view adopted in this paper is that we compute an aliquot sequence until
either we find that it terminates or cycles, or we find that it reaches
some given size. In particular, I pursued every sequence starting
with at most 6 decimal digits to
100 decimal digits (if it did not terminate or cycle before).

The idea of computing aliquot sequences for small starting values $n_0$
is the obvious way to get a feeling for their behaviour,
and hence has been attempted very often. The main problem with this
approach is that for some $n_0$ the values of $s^k(n_0)$ grow rapidly
with $k$; this causes difficulties because all known practical ways
to compute $s(n)$ use the prime factorization of $n$ in an essential
way. Clearly, $s(n)=\sigma(n)-n$, where $\sigma$ denotes
the sum-of-{\it all}-divisors function, which has the advantage over
$s$ of being multiplicative, so it
can be computed using the prime factorization of $n$:
$$\sigma(n)=\prod_{\genfrac{}{}{0pt}{}{p^k\parallel n}{p\text{ prime}}}(1+p+\cdots+p^k),$$
where $p^k\parallel n$ indicates that $p^{k}$ divides $n$ but $p^{k+1}$
does not.

Thus it is no coincidence that similar computations have been performed
over the past 25 years after new factorization algorithms were
developed, and better hardware became much more widely available.
There have been several initiatives following pioneering work of
Wolfgang Creyaufm\"uller \cite{Crey}, and for ongoing progress one
should consult webpages like \cite{Rechenkraft} with contributions by
many individuals.

Despite the extended experience and knowledge gained from
computations such as reported here, it still seems unlikely that
Conjecture (1.1) will be proved or disproved soon; certainly mere computation
will not achieve this. Yet, valuable insight might be obtained.

The current paper grew out of my intermittent attempts over 25 years
to independently perform all necessary computations (at least twice),
and, causing more headaches, to make sure that all confluences
were faithfully recorded.
My main findings are summarized in the table and charts given below.

$$
\vbox{\halign{$\hfil#\hfil$&&\ \hbox to 2.5truecm{$\hfil#\hfil$}\cr
&\text{\it digits} & \text{\it terminating} & \text{\it cycle} & \text{\it open} \cr
& 10 & 735421 & 16204 & 248374 \cr
& 20 & 783786 & 17274 & 198939 \cr
& 30 & 797427 & 17761 & 184811 \cr
& 40 & 800703 & 17834 & 181462 \cr
& 50 & 803317 & 17913 & 178769 \cr
& 60 & 804830 & 17940 & 177229 \cr
& 70 & 805458 & 17985 & 176556 \cr
& 80 & 807843 & 18036 & 174120 \cr
& 90 & 809362 & 18039 & 172598 \cr
& 100& 811555 & 18103 & 170341 \cr
}}$$

\noindent
The table above summarizes what happens if, for starting values up to
$10^6$, we pursue the aliquot sequences up to a size of $d$
decimal digits, with $d$ growing from 10 to 100. As more cycles
and terminating sequences are found, the number of open sequences
declines. In Section \ref{sec:open} a more detailed table is given
for even starting values only.

We try to visualize the rate at which this process takes place in
the pictures below: they plot the number of starting values (below $10^6$)
that terminate or cycle before the given number of digits (on the horizontal
axis) is reached. Note that in both charts the
absolute numbers are plotted vertically, but the scale differs markedly.

\smallskip\noindent
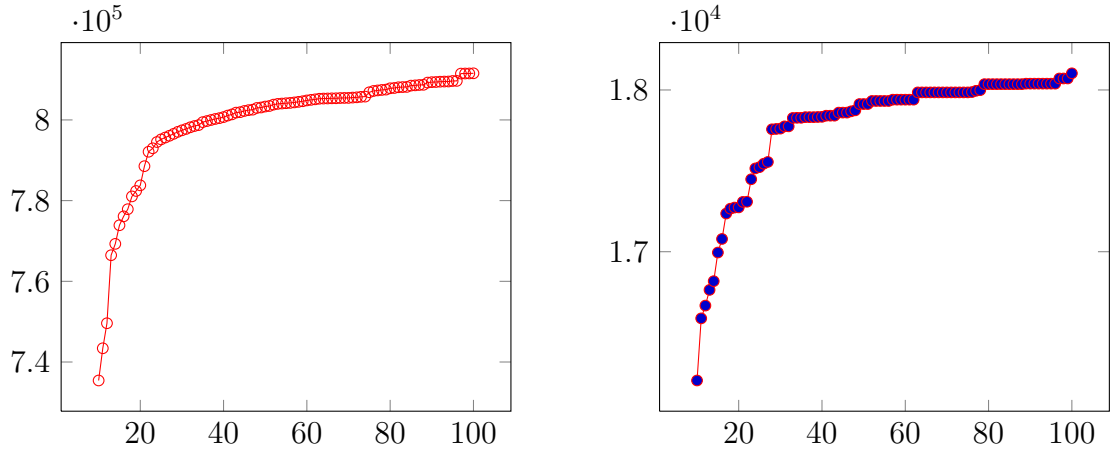
\begin{center}
\hfill
\tikzpicture
\axis[width=7.5cm]
\addplot[color=red,mark=o] coordinates
{
(10, 735421) (11, 743404) (12, 749607) (13, 766467) (14, 769281)
(15, 773897) (16, 776101) (17, 777875) (18, 781046) (19, 782391)
(20, 783786) (21, 788535) (22, 792116) (23, 792972) (24, 794520)
(25, 795181) (26, 795636) (27, 796072) (28, 796547) (29, 797042)
(30, 797427) (31, 797765) (32, 798174) (33, 798482) (34, 798716)
(35, 799455) (36, 799716) (37, 799985) (38, 800238) (39, 800439)
(40, 800703) (41, 801079) (42, 801372) (43, 801813) (44, 801927)
(45, 802226) (46, 802422) (47, 802542) (48, 802966) (49, 803081)
(50, 803317) (51, 803452) (52, 803794) (53, 803975) (54, 804061)
(55, 804139) (56, 804220) (57, 804339) (58, 804480) (59, 804627)
(60, 804830) (61, 804976) (62, 805103) (63, 805252) (64, 805295)
(65, 805318) (66, 805354) (67, 805366) (68, 805428) (69, 805452)
(70, 805458) (71, 805539) (72, 805602) (73, 805723) (74, 805779)
(75, 806856) (76, 807176) (77, 807312) (78, 807432) (79, 807517)
(80, 807843) (81, 807916) (82, 808094) (83, 808141) (84, 808182)
(85, 808505) (86, 808567) (87, 808658) (88, 808703) (89, 809259)
(90, 809362) (91, 809398) (92, 809478) (93, 809511) (94, 809525)
(95, 809652) (96, 809675) (97, 811486) (98, 811506) (99, 811520)
(100, 811555)
};
\endaxis
\endtikzpicture
\hfill
\tikzpicture
\axis[width=7.5cm]
\addplot+[sharp plot,color=red] coordinates
{
(10, 16204) (11, 16588) (12, 16667) (13, 16764) (14, 16818)
(15, 16995) (16, 17079) (17, 17236) (18, 17266) (19, 17273)
(20, 17274) (21, 17308) (22, 17309) (23, 17448) (24, 17515)
(25, 17523) (26, 17544) (27, 17555) (28, 17757) (29, 17759)
(30, 17761) (31, 17775) (32, 17775) (33, 17827) (34, 17828)
(35, 17828) (36, 17832) (37, 17832) (38, 17832) (39, 17833)
(40, 17834) (41, 17841) (42, 17841) (43, 17841) (44, 17860)
(45, 17860) (46, 17860) (47, 17868) (48, 17874) (49, 17913)
(50, 17913) (51, 17913) (52, 17932) (53, 17932) (54, 17932)
(55, 17932) (56, 17932) (57, 17940) (58, 17940) (59, 17940)
(60, 17940) (61, 17940) (62, 17940) (63, 17985) (64, 17985)
(65, 17985) (66, 17985) (67, 17985) (68, 17985) (69, 17985)
(70, 17985) (71, 17985) (72, 17985) (73, 17985) (74, 17985)
(75, 17985) (76, 17987) (77, 17995) (78, 17999) (79, 18036)
(80, 18036) (81, 18036) (82, 18036) (83, 18036) (84, 18036)
(85, 18036) (86, 18036) (87, 18036) (88, 18036) (89, 18039)
(90, 18039) (91, 18039) (92, 18039) (93, 18039) (94, 18039)
(95, 18039) (96, 18039) (97, 18071) (98, 18071) (99, 18071)
(100, 18103)
};
\endaxis
\endtikzpicture
\hfill
\end{center}

The next pair of pictures displays the effect of bounding the starting 
values. In the chart on the left, the bottom graph shows that almost
no starting values less than $10^3$ reach a size of 20 decimal digits,
but for starting values up to $10^4$ around $8\%$ do, a percentage that
grows to more than $13\%$ for starting values up to $10^5$ and $20\%$
up to $10^6$. The corresponding (growing) percentages for terminating
sequences are displayed on the right. The corresponding percentages will almost,
but not exactly, add up to $100\%$, as a small percentage (less than
$2\%$) leads to aliquot cycles.

As a rather naive indication for the truth of the Catalan-Dickson
conjecture, we have also calculated a first order, linear, approximation
($-0.019\cdot x+18.93$) to the percentage of open sequences reaching to more than
20 digits, with starting values up to $10^6$; this the line drawn on
the top left. 

\medskip\noindent
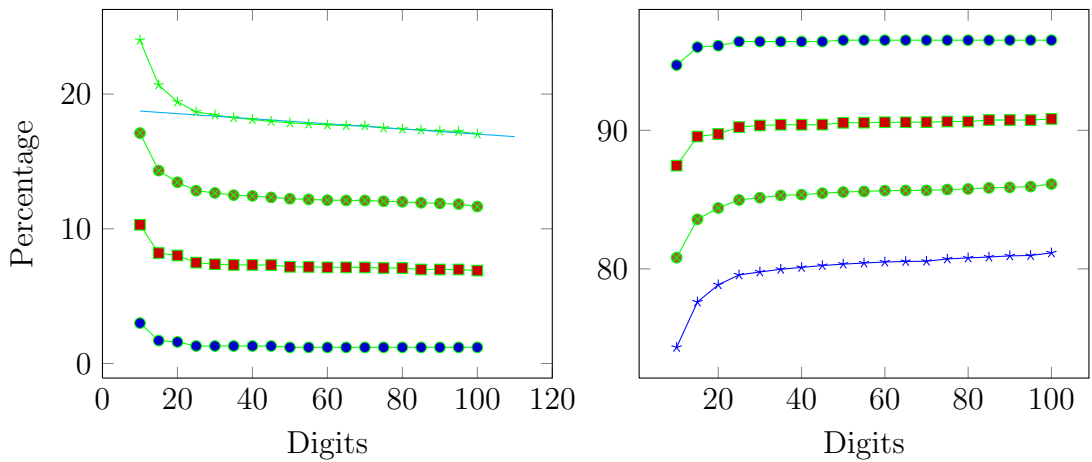
\tikzpicture[domain=10:110]
\axis
[width=7.5cm,xlabel=Digits,ylabel=Percentage]
\addplot+[color=green, sharp plot] coordinates
{
(10, 3.0) (15, 1.7) (20, 1.6) (25, 1.3) (30, 1.3) (35, 1.3) (40, 1.3)
(45, 1.3) (50, 1.2) (55, 1.2) (60, 1.2) (65, 1.2) (70, 1.2) (75, 1.2)
(80, 1.2) (85, 1.2) (90, 1.2) (95, 1.2) (100, 1.2)
};
\addplot+[color=green, sharp plot] coordinates
{
(10, 10.30) (15, 8.20) (20, 8.02) (25, 7.49) (30, 7.37) (35, 7.32) (40, 7.32)
(45, 7.31) (50, 7.18) (55, 7.17) (60, 7.14) (65, 7.14) (70, 7.13) (75, 7.09)
(80, 7.09) (85, 6.98) (90, 6.98) (95, 6.98) (100, 6.90)
};
\addplot+[color=green, sharp plot] coordinates
{
(10, 17.105) (15, 14.320) (20, 13.456) (25, 12.827) (30, 12.663) (35, 12.498)
(40, 12.441) (45, 12.334) (50, 12.230) (55, 12.184) (60, 12.127) (65, 12.109)
(70, 12.105) (75, 12.044) (80, 12.004) (85, 11.929) (90, 11.887) (95, 11.833)
(100, 11.653)
};
\addplot+[color=green, sharp plot] coordinates
{
(10,  24.0007) (15,  20.6819) (20,  19.4156) (25,  18.6819) (30,  18.4459)
(35,  18.2451) (40,  18.1079) (45,  17.9717) (50,  17.8634) (55,  17.7847)
(60,  17.7083) (65,  17.6660) (70,  17.6475) (75,  17.4836) (80,  17.4047)
(85,  17.3396) (90,  17.2562) (95,  17.2285) (100, 17.0322)
};
\addplot [color=cyan, sharp plot] 
{
-0.019*x+18.93
}; 
\endaxis
\endtikzpicture
\tikzpicture[domain=10:110]
\axis
[width=7.5cm,xlabel=Digits] 
\addplot+[color=green, sharp plot] coordinates
{
(10, 94.7) (15, 96.0) (20, 96.1) (25, 96.4) (30, 96.4) (35, 96.4) (40, 96.4)
(45, 96.4) (50, 96.5) (55, 96.5) (60, 96.5) (65, 96.5) (70, 96.5) (75, 96.5)
(80, 96.5) (85, 96.5) (90, 96.5) (95, 96.5) (100, 96.5)
};
\addplot+[color=green, sharp plot] coordinates
{
(10, 87.45) (15, 89.55) (20, 89.73) (25, 90.23) (30, 90.35) (35, 90.40)
(40, 90.40) (45, 90.41) (50, 90.54) (55, 90.55) (60, 90.58) (65, 90.58)
(70, 90.59) (75, 90.63) (80, 90.63) (85, 90.74) (90, 90.74) (95, 90.74)
(100, 90.82)
};
\addplot+[color=green, sharp plot] coordinates
{
(10, 80.813) (15, 83.577) (20, 84.394) (25, 84.976) (30, 85.137) (35, 85.302)
(40, 85.355) (45, 85.462) (50, 85.544) (55, 85.590) (60, 85.647) (65, 85.662)
(70, 85.666) (75, 85.727) (80, 85.767) (85, 85.842) (90, 85.884) (95, 85.938)
(100, 86.118)
};
\addplot+[color=blue,sharp plot] coordinates
{
(10, 74.3404) (15, 77.6101) (20, 78.8535) (25, 79.5636) (30, 79.7765)
(35, 79.9716) (40, 80.1079) (45, 80.2422) (50, 80.3452) (55, 80.4220)
(60, 80.4976) (65, 80.5354) (70, 80.5539) (75, 80.7176) (80, 80.7916)
(85, 80.8567) (90, 80.9398) (95, 80.9675) (100, 81.1574)
};
\endaxis
\endtikzpicture

There is no reason to believe (nor model to support) linear
decay in the long run, but the line does reflect the downward tendency
on the interval between 25 and 100 digits.

Sometimes we find it useful identify aliquot sequences with
the same tails; we say that they {\it merge} at some point.
We call a sequence a {\it main} sequence if it has not merged
with a sequence with a smaller starting value (yet).
The final plot in this section shows data for the number of {\it different}
small aliquot sequences in this sense: the number of aliquot sequences
starting below $10^6$ that exceed $N$ digits for the first time
at different values. At $N=100$ digits there remain 9327 such main
sequences.

\begin{center}
\tikzpicture[domain=10:110]
\axis[ymin=0,ymax=20000]
[width=9cm,xlabel=Digits] 
\addplot+[color=red, sharp plot] coordinates
{
(10, 16782)
(11, 15650)
(12, 14842)
(13, 14160)
(14, 13617)
(15, 13204)
(16, 12866)
(17, 12564)
(18, 12290)
(19, 12071)
(20, 11878)
(21, 11712)
(22, 11564)
(23, 11416)
(24, 11285)
(25, 11180)
(26, 11090)
(27, 11022)
(28, 10947)
(29, 10858)
(30, 10804)
(31, 10732)
(32, 10676)
(33, 10628)
(34, 10571)
(35, 10510)
(36, 10467)
(37, 10428)
(38, 10391)
(39, 10347)
(40, 10304)
(41, 10260)
(42, 10225)
(43, 10197)
(44, 10169)
(45, 10126)
(46, 10095)
(47, 10063)
(48, 10037)
(49, 10011)
(50, 9986)
(51, 9968)
(52, 9937)
(53, 9917)
(54, 9896)
(55, 9874)
(56, 9859)
(57, 9837)
(58, 9815)
(60, 9776)
(65, 9693)
(70, 9633)
(75, 9570)
(80, 9509)
(100, 9327)
};
\endaxis
\endtikzpicture

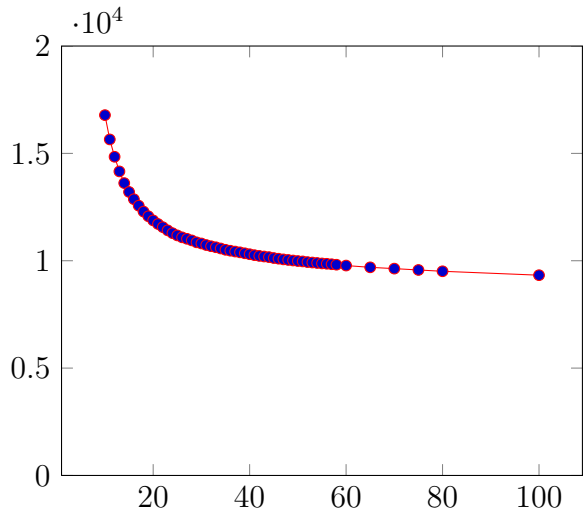
\centerline{\sl Number of remaining main sequences at given number of digits}
\end{center}

\section{Preliminaries}
\noindent
In this section we have collected some known results (with pointers to
the existing literature) as well as some terminology (some standard, some
ad hoc).

Arguments about random integers are not automatically
applicable to heuristics for aliquot sequences due to the fact that certain
factors tend to persist in consecutive values.
The most obvious example of this phenomenon is parity preservation:
$s(n)$ is odd for odd $n$ unless $n$ is an odd square,
$s(n)$ is even for even $n$ unless $n$ is an even square or twice an even square. Guy and Selfridge introduced the notion of driver \cite{GuySelf}.
A {\it driver} of an even integer $n$ is a divisor $2^km$ satisfying
three properties: $2^k\parallel n$;
the odd divisor $m$ is also a divisor of $\sigma(2^k)=2^{k+1}-1$; and,
conversely, $2^{k-1}$ divides $\sigma(m)$.
As soon as $n$ has an additional odd factor (coprime to $v$) besides
the driver, the same driver will also divide $s(n)$. The even perfect numbers
are drivers, and so are only five other integers (2, 24, 120, 672, 523776).
Not only do they tend to persist, but with the exception of $2$,
they also drive the sequence upward, as $s(n)/n$ is $1$ for the
perfect numbers, and $\frac{1}{2}$, $\frac{3}{2}$, $2, 2, 2$ for the other drivers.

More generally, it is possible to prove that arbitrarily long
increasing aliquot sequences exist, a result attributed to H.\,W. Lenstra
(see \cite{Guy}, \cite{Riele}, \cite{Erdos}).

Another heuristic reason to question the truth of the Catalan-Dickson conjecture
was recently refuted in \cite{BK}. We showed that, in the long run,
the growth factor in an aliquot sequence with even starting value
will be less than 1. Besides giving a probabilistic argument
(which does not say anything about counterexamples of `measure 0'),
this is not as persuasive as it may seem, since it assumes that
entries of aliquot sequences behave randomly, which is not true, as
we argued above.

Not only does parity tend to persist in aliquot sequences, the typical
behavior of the two parity classes of aliquot sequences is very different.
There is much stronger tendency for odd $n$ to have
$s(n)<n$. In all odd begin segments, only four cases were encountered
during our computations
where four consecutive odd values were increasing:
\begin{eqnarray*}
&&38745, 41895, 47025, 49695\cr
&&651105, 800415, 1019025, 1070127\cr
&&658665, 792855, 819945, 902295\cr
&&855855, 1240785, 1500975, 1574721.\cr
\end{eqnarray*}
On the other hand, seeing the
factorizations of examples, as in the first quadruple
$$3^3\cdot5\cdot7\cdot41,\quad  3^2\cdot 5\cdot 7^2\cdot 19,\quad
3^2\cdot 5^2\cdot 11\cdot 19,\quad 3\cdot 5\cdot 3313,$$
it is not so difficult to generate longer (and larger!) examples,
such as
$$25399054932615, 37496119518585, 48134213982855, 63887229572985, $$
$$72415060070535, 87397486554105, 101305981941255, 115587206570745, $$
$$133433753777415, 163310053403385, 174881380664583,$$
in the vein of the result of Lenstra, but such examples
did not occur in our sequences yet.

We say that sequence $s$ {\it merges} with sequence $t$
(at value $x$) if $s$ and $t$ have $x$ as first common value,
$t$ has a smaller starting value than $s$,
and the common value occurs before $s$ reaches its maximum.
In this case $t$ will be the {\it main} sequence (unless
it merges with a `smaller' sequence again).
From $x$ on, $s$ and $t$ will coincide of course.

We should point out again that the notion of being a main sequence
is time dependent: sequences may merge beyond
the point to which we have as yet computed them.

Since all of our sequences are finite (except, possibly, for a repeating
cycle at the end), we can speak of
the {\it height} of a sequence: this is essentially the logarithm
of its maximal value; sometimes we measure this in number of decimal
digits, sometimes in number of bits. The {\it volume} will be the
sum of the number of digits of the entries of the sequence, without
rounding first: so $\vol(s) = \sum_{x\in s} \log_{10} x$.

\section{Odd cases}
\noindent
We first consider the 500000 odd starting values,
as they usually lead to termination quickly.
In fact, 494088 odd starting
values terminated; 5119 sequences with odd starting values lead
to a cycle (see next section) and 793 remain open, after merging
with a sequence with an even starting value.

As we saw above, parity is not always maintained. Therefore we need to
distinguish in our 500000 odd starting values between
aliquot sequences consisting {\it only}
of odd integers, and those containing even values as well.

For 440239 odd starting values, an all-odd sequence ensues;
the remaining 59761 change over to even, after hitting an odd square.
Of the 59761 odd starting values that change over, 12674 do so after hitting
$3^2$. Only the odd squares less than a million did occur.

Of the all-odd sequences, 208 end in an odd cycle.

Of the 59761 odd-to-even starting values, 5119 lead to a cycle
and 793 merge with an even sequence reaching 100 digits.
Of the 54057 terminating odd-to-even starting values,
17 take more than 1000 steps before terminating:
11 of them merge after a couple of steps with the 94-digit maximum length
1602 sequence 16302, and 6 of them merge after a couple of steps
with the 76-digit-maximum length 1740 sequence 31962.

Here are the numbers in summary:

$$
\vbox{\halign{$\hfil#\hfil$&&\ $\hfil#\hfil$\cr
&\text{\it parity} & \text{\it terminating} & \text{\it cycle} & \text{\it open} \cr
& \text{\rm odd} & 494088 & 5119 & 793 \cr
& \text{\rm all-odd} & 440031 & 208 & 0 \cr
& \text{\rm even} & 317467 & 12984 & 169548 \cr
& \text{\rm all} & 811555 & 18103 & 170341 \cr
}}$$

\section{Terminators}
\noindent
Of the 999999 starting values, 811555 terminated without reaching
a 100-digit value. 
\subsection{All-odd terminators}
Among the 440031 all-odd terminators, 78497 terminate after 1 step. This
reflects that there are 78497 odd primes less than a million. The
longest all-odd terminator has length 23:
\begin{eqnarray*}
\section{Open}\label{sec:open}
\noindent
The following table breaks up the range of starting values into
ten sub-intervals from $k\cdot 10^5$ to $(k+1)\cdot 10^5-1$,
for $k=0, 1, \ldots, 9$, and for those the number of even starting
values reaching $d$ decimal digits is given. 
Note that (as 0 is not included as starting value) the first column
concerns 49999 starting values, and the other columns 50000.
The final column is the sum of the first 10 columns, and counts how
many of the 499999 even starting values less than $10^6$ reach
$d$ decimal digits (that is, a value of at least $10^{99}$).
\noindent
the 461214 sequence merges with the open 4788 sequence after 6467 steps
(after reaching a 88 digit local maximum). To complicate the situation,
it first merges with the 314718 sequence
(461214:5=314718:4=1372410) which in turn merges with the 4788 sequence
(on its way picking up 14 more sequences that have the same local maximum.
The longest of these, the 461214 and 580110 sequences, reach 100 digits (with 4788)
after 8599 steps. The next longest pre-merger example is a group of 4
sequences merging with the open 1920 sequence after 4656 steps and
a 76 digit maximum.
The tables below list all cycles that occur, with their popularity.
The second column lists the number of starting values ending in the
cycle listed in the first column, with (in parentheses) the number
of {\it main} sequences among these.
The third column lists the number of even starting values among those
of the second column. In the fourth column is shown how often each
of the entries of the cycle is first hit by some sequence. Thus, for example,
the entry $2\ /\ 9$ in the row for the amicable pair $[220, 284]$
reflects that besides the starting values 220 and 284 only 9 other
sequences up to $10^6$ lead to this cycle (8 of them with odd starting
value according to column 3) and only one of those will hit 220 first.

\begin{small}
$$\vbox{\halign{{$#\hfil$} : 
&\qquad \hfil$#$ 
&\qquad \hfil$#$
&\qquad \hfil$#$
&\qquad \hfil$#$\hfil
&\qquad \hfil$#$\hfil
&\qquad \hfil$#$\hfil
&\qquad \hfil$#$\hfil
&\qquad \hfil$#$\hfil\cr
\text{cycle} & \text{\#} & \text{(\#main)} & \text{even} & \text{entry} \cr
&  & & & \cr
[\ 6\ ] &  5395 & (5132) & 579 & 5395      \cr
[\ 28\ ] &  1& (1) & 1 & 1     \cr
[\ 496\ ] &  13 & (11) & 12 & 13     \cr
[\ 8128\ ] &  1408 & (460) & 1408 & 1408     \cr
[ 1264460, 1547860,  & & & & \cr
\ \ \ \ 1727636, 1305184 ] &  13& (2) & 13 & 13\vert 0\vert 0\vert 0  \cr
[ 2115324, 3317740,  &  1& (1) & 1 & 1\vert 0\vert 0\vert 0  \cr
\ \ \ \ 3649556, 2797612 ] &  1& (1) & 1 & 1\vert 0\vert 0\vert 0  \cr
[ 12496, 14288,  &  & & & \cr
\ \ \ \  15472, 14536, 14264 ] &  150& (109) & 150 & 72\vert 2\vert 1\vert 74\vert 1 \cr
C_{28}
& 741 & (131) & 741 & 8\vert 1\vert 3\vert 3\vert 1\vert 6\vert 1\vert 2\cr 
& & & & 33\vert 1\vert5\vert 1\vert 2\vert 19\vert 15 \cr 
& & & & 1\vert 157\vert 1\vert 1\vert 1\vert 3\vert 5 \cr 
& & & & 1\vert 35\vert 1\vert 49\vert 269\vert 123 \cr 
\text{\rm total} &  18103 & (11056) & 12984 \cr
}}$$

\smallskip\noindent
$$\vbox{\halign{$#\hfil$ : 
&\qquad \hfil$#$ 
&\qquad \hfil$#$
&\qquad \hfil$#$\hfil
&\qquad \hfil$#$\hfil
&\qquad \hfil$#$\hfil
&\qquad \hfil$#$\hfil
&\qquad \hfil$#$\hfil\cr
\text{cycle} & \text{\#}\hfill\text{(\#non-merging)} & \text{even} & \text{entry} \cr
[ 220, 284 ] &  11 \hfill (10) & 3 & 1\vert 10    \cr
[ 1184, 1210 ] &  7564 \hfill (3841) & 7561 & 3599\vert 3965    \cr
[ 2620, 2924 ] &  1153 \hfill (533) & 1152 & 9\vert 1144    \cr
[ 5020, 5564 ] &  50 \hfill (44) & 24 & 1\vert 49    \cr
[ 6232, 6368 ] &  27 \hfill (26) & 10 & 26\vert 1    \cr
[ 10744, 10856 ] &  249 \hfill (125) & 249 & 1\vert 248    \cr
[ 12285, 14595 ] &  106 \hfill (104) & 0 &  56\vert 50    \cr
[ 17296, 18416 ] &  202 \hfill (100) & 202 & 200\vert 2    \cr
[ 63020, 76084 ] &  9 \hfill (2) & 9 & 1\vert 8    \cr
[ 66928, 66992 ] &  6 \hfill (5) & 6 & 5\vert 1    \cr
[ 67095, 71145 ] &  47 \hfill (45) & 0 & 43\vert 4    \cr
[ 69615, 87633 ] &  39 \hfill (36) & 0 & 21\vert 18    \cr
[ 79750, 88730 ] &  342 \hfill (102) & 303 & 306\vert 36    \cr
[ 100485, 124155 ] &  4 \hfill (3) & 0 & 2\vert 2    \cr
[ 122265, 139815 ] &  3 \hfill (2) & 0 & 2\vert 1    \cr
[ 122368, 123152 ] &  3 \hfill (2) & 3 & 2\vert 1    \cr
[ 141664, 153176 ] &  10 \hfill (6) & 10 & 1\vert 9    \cr
[ 142310, 168730 ] &  5 \hfill (4) & 5 & 1\vert 4    \cr
[ 171856, 176336 ] &  23\hfill (17) & 23 & 8\vert 15    \cr
[ 176272, 180848 ] &  17\hfill (7) & 17 & 16\vert 1    \cr
[ 185368, 203432 ] &  106\hfill (56) & 106 & 102\vert 4    \cr
[ 196724, 202444 ] &  25\hfill (19) & 25 & 6\vert 19    \cr
[ 280540, 365084 ] &  121\hfill (41) & 121 & 120\vert 1    \cr
[ 308620, 389924 ] &  6\hfill (5) & 6 & 5\vert 1    \cr
[ 319550, 430402 ] &  17\hfill (8) & 17 & 15\vert 2    \cr
[ 356408, 399592 ] &  2\hfill (1) & 2 & 1\vert 1    \cr
[ 437456, 455344 ] &  12\hfill (6) & 12 & 2\vert 10    \cr
[ 469028, 486178 ] &  34\hfill (10) & 34 & 30\vert 4    \cr
[ 503056, 514736 ] &  9\hfill (5) & 9 & 8\vert 1    \cr
[ 522405, 525915 ] &  6\hfill (5) & 0 & 4\vert 2    \cr
[ 600392, 669688 ] &  3\hfill (2) & 3 & 1\vert 2    \cr
[ 609928, 686072 ] &  3\hfill (1) & 3 & 2\vert 1    \cr
[ 624184, 691256 ] &  5\hfill (1) & 5 & 3\vert 2    \cr
[ 635624, 712216 ] &  39\hfill (10) & 39 & 31\vert 8    \cr
[ 643336, 652664 ] &  2\hfill (1) & 2 & 1\vert 1    \cr
[ 667964, 783556 ] &  7\hfill (5) & 7 & 4\vert 3    \cr
[ 726104, 796696 ] &  4\hfill (3) & 4 & 3\vert 1    \cr
[ 802725, 863835 ] &  2\hfill (1) & 0 & 1\vert 1    \cr
[ 879712, 901424 ] &  35\hfill (4) & 35 & 15\vert 20    \cr
[ 898216, 980984 ] &  9\hfill (1) & 9 & 8\vert 1    \cr
[ 947835, 1125765 ] &  1\hfill (1) & 0 & 1\vert 0    \cr
[ 998104, 1043096 ] &  2\hfill (1) & 2 & 2\vert 0    \cr
[ 1077890, 1099390 ] &  19\hfill (1) & 19 & 19\vert 0    \cr
[ 2723792, 2874064 ] &  13\hfill (3) & 13 & 9\vert 4    \cr
[ 4238984, 4314616 ] &  16\hfill (1) & 16 & 0\vert 16    \cr
[ 4532710, 6135962 ] &  6\hfill (1) & 6 & 6\vert 0    \cr
[ 5459176, 5495264 ] &  6\hfill (1) & 6 & 6\vert 0    \cr
[ 438452624, 445419376 ] &  1\hfill (1) & 1 & 1\vert 0    \cr
}}$$
\end{small}

\end{document}